\documentclass[10pt,twoside]{article}
\usepackage{graphicx}
\usepackage{amsmath}
\usepackage{latexsym, amsfonts,amssymb}
\usepackage{Latex-document}

\title{\bf  Automorphic \boldmath$L$-Functions  \vskip -2mm
and Functoriality\vskip 6mm}
\author{Freydoon Shahidi\vspace*{-0.5cm}\thanks{Purdue University, Department
of Mathematics, West Lefayette, Indiana 47907, USA. E-mail:
shahidi@math.purdue.edu}}
\date{\vspace{-8mm}}

\begin{document}

\maketitle

\thispagestyle{first} \setcounter{page}{655}

\begin{abstract}

\vskip 3mm

This is a report on the global aspects of the Langlands-Shahidi method which in conjunction with converse theorems
of Cogdell and Piatetski-Shapiro has recently been instrumental in establishing a significant number of new and
surprising cases of Langlands Functoriality Conjecture over  number fields. They have led to striking new
estimates  towards Ramanujan and Selberg conjectures.

\vskip 4.5mm


\end{abstract}

\vskip 12mm

\section{Preliminaries} \label{section 1}\setzero

\vskip-5mm \hspace{5mm}

Let $F$ be a number field. For each place $v$ of $F$, let $F_v$ be
its completion at $v$. Assume $v$ is a finite place and let $O_v$
denote the ring of integers of $F_v$. Denote by $P_v$ its maximal
ideal and fix a uniformizing parameter $\varpi_v$ generating
$P_v$. Let $[O_v:P_v]=q_v$ and fix and absolute value $|\quad |_v$
for which $|\varpi_v|_v=q_v^{-1}$.

Let $\bf{G}$ be a quasisplit connected reductive algebraic group
over $F$. Fix an $F$-Borel subgroup $\bf{B=TU}$, where $\bf{T}$ is
a maximal torus of $\bf{B}$ and $\bf{U}$ is its unipotent radical.
Let $\bf{A}_0\subset\bf{T}$ be the maximal split subtorus of
$\bf{T}$. Throughout this article, $\bf{P}$ is a maximal parabolic
subgroup of $\bf{G}$, defined over $F$, with a Levi decomposition
$\bf{P}=\bf{MN}$, where $\bf{M}$ is a Levi subgroup of $\bf{P}$
and $\bf{N}$ is its unipotent radical. We will assume $\bf{P}$ is
standard in the sense that $\bf{N}\subset\bf{U} $. We fix $\bf{M}$
by assuming $\bf{T}\subset\bf{M}$. We finally use $W$ to denote
the Weyl group of $\bf{A}_0$ in $\bf{G}$.

Let $\mathbb{A}_F$ denote the ring of adeles of $F$ and for every
algebraic group $\bf{H}$ over $F$, let $H={\bf{H}}(\mathbb{A}_F)$.
Considering $\bf{H}$ as a group over each $F_v$, we then set
$H_v={\bf{H}}(F_v)$.

Let $\bf{A}$ denote the split component of $\bf{M}$, i.e., the
maximal split subtorus of the connected component of the center of
$\bf{M}$. For every group $\bf{H}$ defined over $F$, let
$X({\bf{H}})_F$ be the group of $F$-rational characters of
$\bf{H}$. We set ${\frak a}
=\mathrm{Hom}(X({\bf{M}})_F,\mathbb{R})$. Then ${\frak
a}^*=X({\bf{M}})_{F\otimes\mathbb{Z}}\mathbb{R}=X({\bf{A}})_{F\otimes\mathbb{Z}}
\mathbb{R}$ and ${\frak a}_{\mathbb{C}}^{\ast}={\frak
a}^{\ast}\otimes_{\mathbb{R}}\mathbb{C}$ is the complex dual of
${\frak a}$.

When $\bf{G}$ is unramified over a place $v$, we let
$K_v={\bf{G}}({\bf{O}}_v)$. Otherwise, we shall fix a special
maximal compact subgroup $K_v\subset G_u$ for which
$G_v=P_vK_v=B_vK_v$. Let $K={\otimes}_vK_v$ Then $G=PK=BK$. Let
$K_M=K\cap M$.

For each $v$, the embedding $X({\bf{M}})_F\hookrightarrow X({\bf{M}})_{F_v}$ induces a map
$$a_v=\mathrm{Hom}(X({\bf{M}})_{F_v},\mathbb{R})\rightarrow {\frak a}.$$
There exists a homomorphism
$H_M:M\rightarrow {\frak a}$ defined by
$$\exp\langle\chi,H_M(m)\rangle=\prod_v|\chi(m_v)|_v$$
for every $\chi\in X({\bf{M}})_F$ and $m=(m_v)$. We extend $H_M$ to $H_P$ on $G$ by making it trivial on $N$ and
$K$.

Let $\alpha$ denote the unique simple root of $\bf{A}$ in
$\bf{N}$. It can be identified by a unique simple root of
$\bf{A}_0$ in $\bf{U}$. If $\rho_{\bf{P}}$ is half the sum of
$F$-roots in $\bf{N}$, we set
$\widetilde{\alpha}=\langle\rho_{\bf{P}},\alpha\rangle^{-1}\rho_{\bf{P}}\in
a^*$, where for each pair of non-restricted roots $\alpha$ and
$\beta$ of $\bf{T}$,
$\langle\alpha,\beta\rangle=2(\alpha,\beta)/(\beta,\beta)$ is the
Killing form.

Given a connected reductive algebraic group $\bf{H}$ over $F$,
let ${}^LH$ be its $L$-group. Considering $\bf{H}$ as a group over
$F_v$, we then denote by ${}^LH_v$ its $L$-group over $F_v$. Let
${}^LH^0={}^LH_v^0$ be the corresponding connected component of
1. We then have a natural homomorphism from ${}^LH_v$ into
${}^LH$. We let $\eta_v:{}^LM_v\rightarrow{}^LM$ be this map for
$M$ (cf. [4]).

Let ${}^LN$ be the $L$-group of $\bf{N}$ defined naturally in [4].
Let ${}^L{\frak n}$ be its ( complex ) Lie algebra, and let $r$
denote the adjoint action
 of ${}^LM$ on ${}^L{\frak n}$. Decompose $r=\bigoplus\limits_{i=1}^m r_i$ to its irreducible
 subrepresentations, indexed according to the values $\langle\widetilde{\alpha},\beta\rangle=i$
 as $\beta$ ranges among the positive roots of $\bf{T}$. More precisely,
 $X_{\beta\vee}\in{}^L{\frak n}$ lies in the space of $r_i$ if and only if
 $\langle\widetilde{\alpha},\beta\rangle=i$. Here $X_{\beta\vee}$ is a root vector
attached to the coroot $\beta^{\vee}$, considered as a root of the
$L$-group. The integer $m$ is equal to the nilpotence class of
${}^L{\frak n}$. We let $r_{i,v}=r_i\cdot\eta_v$ for each $i$
(cf. [34,40,41]).

If $\Delta$ denotes the set of simple roots of ${\bf{A}}_0$ in
$\bf{U}$, we use $\theta\subset{\Delta}$ to denote the subset
generating M. Then $\Delta=\theta\cup \{\alpha \}$. There exists a
unique element $\widetilde{w}_0\in W$ such that
$\widetilde{w}_0(\theta)\subset\Delta$, while
$\widetilde{w}_0(\alpha)<0$. We will always choose a
representative $w_0$ for $\widetilde{w}_0$ in ${\bf{G}}(F)$ and
use $w_0$ to denote each of its components.

\section{Eisenstein series and \boldmath$L$-functions} \label{section 2}\setzero

\vskip-5mm \hspace{5mm}

Let $\pi=\otimes_v\pi_v$ be a cusp form on $M$. Given a
$K_M$-finite function $\varphi$ in the space of $\pi$, we extend
$\varphi$ to a function $\widetilde{\varphi}$ on $G$ as in Section
2 of [39] as well as in [17], and for $s\in\mathbb{C}$, set
\begin{equation} \label{2.1}
  \phi_s(g)=\widetilde{\varphi}(g)\exp \langle
  s\widetilde{\alpha}+\rho_{\bf{P}},H_P(g) \rangle.
\end{equation}
The corresponding Eisenstein series is then defined by
\begin{equation} \label{2.2}
  E(s,\phi_s,g,P)=\sum_{\gamma\in{\bf{P}}(F) \backslash
  {\bf{G}}(F)}\phi_s(\gamma g)
\end{equation}
(cf. [17,33,34,35]).

Let $I(s,\pi)=\otimes_vI(s,\pi_v)$ be the representation
parabolically induced from $\pi\otimes\exp \langle
s\widetilde{\alpha},H_p(\, ) \rangle $.

Let $\bf{M}'$ be the Levi subgroup of $\bf{G}$ generated by
$\widetilde{w}(\theta)$. There exists a parabolic subgroup
${\bf{P}}'\supset{\bf{B}}$ which has $\bf{M}'$ as a Levi factor.
Let $\bf{N}'$ be its unipotent radical. Given $f$ in the space of
$I(s,\pi)$ and $\mathrm{Re}(s)>>0$, define the global intertwining
operator $M(s,\pi)$ by
\begin{equation} \label{2.3}
  M(s,\pi)f(g)=\int_{N'}f(w_0^{-1}n'g)dn'\quad (g\in G).
\end{equation}
Observe that if $f=\otimes_vf_v$, then for
almost all $v$, $f_v$ is the unique $K_v$-fixed functions normalized by $f_v(e_v)=1$. Finally, if at each $v$ we
define a local intertwining operator by
\begin{equation} \label{2.4}
  A(s,\pi_v,w_0)f_v(g)=\int_{N'_v}f_v(w_0^{-1}n'g)dn',
\end{equation}
then
\begin{equation} \label{2.5}
  M(s,\pi)=\otimes_vA(s,\pi_v,w_0).
\end{equation}

 It follows form the general theory of Eisenstein
series that the poles of\linebreak $E(s,\widetilde{\varphi},g,P)$,
as $\widetilde{\varphi}$ and $g$ vary, are the same as those of
$M(s,\pi)$, and for $\mathrm{Re}(s)\geq 0$, they are all simple
and finite in number, with none on the line $\mathrm{Re}(s)=0$
(cf. [17,33,35]).

By construction each $\phi_s$ belongs to the space of $I(s,\pi)$.
Consequently, one can consider $M(s,\pi)\phi_s$ which is a member
of $I(-s,w_0(\pi))$. The Eisenstein series
$E(s,\widetilde{\varphi},g,P)$ then satisfies the functional
equation
\begin{equation} \label{2.6}
  E(s,\phi_s,g,P)=E(-s,M(s,\pi)\phi_s,g,P').
\end{equation}

Suppose that $\bf{G}$ splits over $L$, where $L$ is a finite
Galois extension of $F$. For every unramified $v$, there exists a
unique Frobenius conjugacy class in Gal($L_w/F_v$), $w|v$ which we
denote by $\tau_v$. Moreover, if $v$ is such that $\pi_v$ and
$\bf{G}$ are both unramified, then there exists and ${}^LM$
semisimple conjugacy class in ${}^LM^0\rtimes\tau_v$ which
determines $\pi_v$ uniquely ([4,40]). We may identify, as we in
fact do, this conjugacy class with an element $A_v\in{}^LT^0$
which may be assumed to be fixed by $\tau_v$ (cf. \S6.3 and 6.5 of
[4]). The local Langlands $L$-function defined by $\pi_v$ and
$r_v,r_v=r\cdot\eta_v$, where $r$ is a complex analytic
representation of ${}^LM$, is then defined to be (cf. [4,34,40]),
\begin{equation} \label{2.7}
  L(s,\pi_v,r_v)=\det(I-r_v(A_v\rtimes\tau_v)q^{-1}_v)^{-1}.
\end{equation}

Let $S$ be a finite set of places of $F$, including all the
archimedean ones, such that for every $v\notin S$, $\pi_v$ and
$\bf{G}$ are both unramified. Set
\begin{equation} \label{2.8}
  L_S(s,\pi,r)=\prod_{v\notin S}L(s,\pi_v,r_v).
\end{equation}
The main result of [34, also see 40] is that
\begin{eqnarray}
  M(s,\pi)f&=&\otimes_{v\in S}A(s,\pi_v,w_0)f_v\otimes\otimes_{v\notin S}\widetilde{f}_v \nonumber \\
  & & \times \prod_{i=1}^mL_S(is,\pi,\widetilde{r}_i)/L_S(1+is,\pi,\widetilde{r}_i), \label{2.9}
\end{eqnarray}
where $f=\otimes_vf_v$ is such that for each $v\notin S$, $f_v$
is the unique $K_v$-fixed function in $I(s,\pi_v)$ normalized by
$f_v(e_v)=1$ and for each $i$, $\widetilde{r}_i$ denotes the
contragredient of $r_i, i=1, \cdots, m$, the irreducible
components of the adjoint action of ${}^LM$ or ${}^LN$. Here
$\widetilde{f}_v$ is the $K_v$-fixed function in the space of
$I(-s,w_0(\pi_v))$, normalized the same way. Moreover $f_v$ and
$\widetilde{f}_v$ are identified as elements in spherical
principal series.

\section{Generic representations and the non-constant term} \label{section 3}\setzero

\vskip-5mm \hspace{5mm}

Suppose $F$ is a field, either local or global, and $\bf{G}$ is as
before, with a Borel subgroup $\bf{B=TU}$ over $F$. Fix an
$F$-splitting \{$X_{\alpha'}$\}, i.e., a collection of root
vectors as $\alpha'$ ranges over simple roots of $\bf{T}$ in
$\bf{U}$ which is invariant under the action of
$\Gamma_F=\mathrm{Gal}(\overline{F}/F)$. This then determines a
map $\phi$ form $\bf{U}$ to $\Pi
\mathbb{G}_a,\varphi(u)=(x_{\alpha'})_{\alpha'}$, where
$x_{\alpha'}$ is the $\alpha'$-coordinate of $u$ with respect to
$\{X_{\alpha'}\}$. Let $\{\kappa_{\alpha'}\}$ be a collection of
elements in $\overline{F}^*$ such that
$\sigma(\kappa_{\alpha'})=\kappa_{\sigma\alpha'}$ for every
$\sigma\in \Gamma_F$. Set
$f(u)=\sum\limits_{\alpha'}\kappa_{\alpha'}x_{\alpha'}$. Observe
that $f$ is $F$-rational. If $F$ is global, we extend $f$ to a
map on ${\bf{U}}({\mathbb{A}}_F)$. Let $\psi_F$ be a non-trivial
character of $F$ ($F\setminus\mathbb{A}_F$ if $F$ is global). A
character $\chi$ of
${\bf{U}}(F)({\bf{U}}(F)\setminus{\bf{U}}(\mathbb{A}_F)\ \mbox
{if}\ F\ \mbox {is global})$ is called $non-degenerate$ or
$generic$ if
$\chi(u)=\varphi(f(u)),u\in{\bf{U}}(F)(u\in{\bf{U}}(F)\setminus{\bf{U}}(\mathbb{A}_F)$
if $F$ is global).

We now continue to assume $F$ is a number field. Let
$\chi=\otimes_v\chi_v$ be a generic character of
${\bf{U}}(F)\setminus U$.

Let ${\bf{U}}^0={\bf{U}}\cap{\bf{M}}$ and let $\chi$ also denote
the restriction of $\chi$ to $U^0$. Choose a function $\varphi$ in
the space of $\pi=\otimes_v\pi_v$, a cuspidal representation of
$M$, and ${\bf{U}}^0(F)\setminus U^0$ being compact, set
\begin{equation} \label{3.1}
  W_{\varphi}(m)=\int_{{\bf{U}}^0(F)\setminus U^0}\varphi(um)\overline{\chi(u)}du.
\end{equation}
We shall say $\pi$ is (globally) $\chi$-generic if
$W_{\varphi}\neq 0$ for some $\varphi$. The representation $\pi$
is (globally) generic if it is $\chi$-generic with respect to
some generic $\chi$. Then each $\pi_v$ will be $\chi_v$-generic
in the sense that there exists a non-zero Whittaker functional
$\lambda_v$ i.e., a continuous (in the semi-norm topology if
$v=\infty$) functional satisfying
$\langle\pi_v(u)x,\lambda_v\rangle=\chi_v(u)\langle
 x,\lambda_v\rangle,x\in \mathcal {H}(\pi_v),u\in U_v^0$. Choosing
$\varphi$ appropriately, i.e., if
$\varphi=\otimes_v\varphi_v,\varphi_v\in \mathcal {H}(\pi_v)$,
then
$W_\varphi(m)=\prod_v\langle\pi_v(m_v)\varphi_v,\lambda_v\rangle$,
for $m=(m_v)$.

Given $f_v\in V(s,\pi_v)$, the space of $I(s,\pi_v)$, define
\begin{equation} \label{3.2}
  \lambda_{\lambda_v}(s,\pi_v)(f_v)=\int_{N'_v}\langle f_v(w_0^{-1}n'),\lambda_v
  \rangle\overline{\chi(n')}dn',
\end{equation}
a canonical Whittaker functional for $I(s,\pi_v)$. Changing the
splitting we now assume $\kappa_{\alpha'}=1$. It now follows from
Rodier's theorem that there exists a complex function (of\,$s$),
$C_{\chi_v}(s,\pi_v)$, depending on $\pi_v,\chi_v$ and $w_0$ such
that (cf. [41,42,43])
\begin{equation} \label{3.3}
  \lambda_{\chi_v}(s,\pi_v)=C_{\chi_v}(s,\pi_v)\lambda_{\chi_v}(-s,w_0(\pi_v))\cdot
  A(s,\pi_v,w_0).
\end{equation}
This is what we call the \it{Local Coefficient} \rm attached to $s,\pi_v,\chi_v$ and $w_0$. The choice of $w_0$ is
now specified by our fixed splitting as in [43].

Finally, if
\begin{equation} \label{3.4}
  E_\chi(s,\phi_s,g,P)=\int_{{\bf{U}}(F)\setminus U}E(s,\phi_s,ug,P)\overline{\chi(u)}du
\end{equation}
is the $\chi$-nonconstant term of the Eisenstein series, then ([7,41,42])
\begin{equation} \label{3.5}
  E_\chi(s,\phi_s,e,P)=\prod_{v\in S}W_v(e)\prod_{i=1}^mL_S(1+is,\pi,\widetilde{r}_i)^{-1},
\end{equation}
where now $S$ is assumed to have the property that if $v\notin S$, then $\chi_v$ is also unramified.

Applying Definition (3.4) to both sides of (2.6), using (3.5) now
implies the \it{crude functional equation} \rm([40,41])
\begin{equation} \label{3.6}
  \prod_{i=1}^mL_S(is,\pi,r_i)=\prod_{v\in S}C_{\overline{\chi}_v}(s,\widetilde{\pi}_v)\prod_{i=1}^m
  L_S(1-is,\pi,\widetilde{r}_i).
\end{equation}

\section{The main induction, functional equations and multiplicativity} \label{section 4}\setzero

\vskip-5mm \hspace{5mm}

To prove the functional equation for each $r_i$ with precise root
numbers and $L$-function, we use (cf. [42]):

{\bf Proposition 4.1.} \it Given $1<i\leq m$, there exists a
quasisplit guoup ${\bf{G}}_i$ over F, a maximal F-parabolic
subgroup ${\bf{P}}_i={\bf{M}}_i{\bf{N}}_i$, both unramified for
every $v\notin S$, and a cuspidal automorphic form $\pi'$ of
$M_i={\bf{M}}_i({\mathbb{A}}_F)$, unramified for every $v\notin
S$, such that if the adjoint action $r'$ of ${}^LM_i$ on
${}^L{\frak n}_i$ decomposes as
$r'=\bigoplus\limits_{j=1}^{m'}r'_j$, then
$$L_S(s,\pi,r_i)=L_S(s,\pi',r'_1).$$
Moreover $m'<m$. \rm

{\bf Remark 4.2.} As was observed by Arthur [1], each
${\bf{M}}_i$ can be taken equal to $\bf{M}$ and $\pi'=\pi$. In
fact each ${\bf{G}}_i$ can be taken to be an endoscopic group for
${\bf{G}}$, sharing $\bf{M}$ as a Levi subgroup. We shall record
this as

{\bf Proposition 4.3.} \it Given $i$, $1<i\leq m$, there exist a quasisplit connected reductive $F$-group with
$\bf{M}$ as a Levi subgroup and $m'<m$ for which $r'_1=r_i$.\rm

Using this induction and local-global arguments (cf. Proposition
5.1 of [42]), it was proved in [42] that

{\bf Theorem 4.4.}  $($Theorems 3.5 and 7.7 of [42]$)$ {\it $a)$
For each $i$, $1\leq i\leq m$, and each $v$, there exist a local
$L$-function $L(s,\pi_v,r_{i,v})$, which is the inverse of a
polynomial in $q_v^{-s}$ whose constant term is 1, if $v<\infty$,
and is the Artin $L$-function attached to $r_i\cdot\varphi'_v$,
where $\varphi'_v:W'_{F_v}\rightarrow{}^LM_v$ is the homomorphism
of the Deligne-Weil group into ${}^LM_v$ parametrizing $\pi_v$,
if either $v=\infty$ or $\pi_v$ has an Iwahori-fixed vector; and
a root number $\varepsilon(s,\pi_v,r_{i,v},\varphi_v)$ satisfying
the same provisions, such that if
\begin{equation} \label{4.4.1}
  L(s,\pi,r_i)=\prod L(s,\pi_v,r_{i,v})
\end{equation}
and
\begin{equation} \label{4.4.2}
  \varepsilon(s,\pi,r_i)=\prod_v \varepsilon(s,\pi_v,r_{i,v},\psi_v),
\end{equation}
then
\begin{equation} \label{4.4.3}
  L(s,\pi,r_i)=\varepsilon(s,\pi,r_i)L(1-s,\pi,\widetilde{r}_i).
\end{equation}

$b)$ Let
\begin{equation} \label{4.4.4}
  \gamma(s,\pi_v,r_{i,v},\psi_v)=\varepsilon(s,\pi_v,r_{i,v},\psi_v)L(1-s,\pi_v,
  \widetilde{r}_{i,v})/L(s,\pi_v,r_{i,v}).
\end{equation}
Then each $\gamma(s,\pi_v,r_{iv},\psi_v)$ is multiplicative in the sense of equation (3.13) in Theorem 3.5 of
[42]. (See below.) If $\pi_v$ is tempered, then $\gamma(s,\pi_v,r_{i,v},\psi_v)$ determines the corresponding root
number and L-function uniquely and in fact that is how they are defined. Suppose $\pi_v$ is non-tempered, then
each $L(s,\pi_v,r_{i,v})$ is determined by means of the analytic continuation of its quasi-tempered Langlands
parameter and multiplicativity of corresponding $\gamma$-functions. More precisely, if $\sigma_v$ is the
quasitempered Langlands parameter that gives $\pi_v$ as a subrepresentation, then
\begin{equation} \label{4.4.5}
  L(s,\pi_v,r_{i,v})=\prod_{j\in S_i}L(s,\overline{w}_j(\sigma_v),r'_{i(j),v}),
\end{equation}
where the notation is as in part 3) of Theorem 3.5 of [42],
provided that every L-function on the right hand side is
holomorphic for $\mathrm{Re}(s)>0$, whenever $\sigma_v$ is
(unitary) tempered (Conjecture 7.1 of [42], proved in many cases
[3.6.42]). The set $S_i, \overline{w}_j$ and $r'_{i(j)}$ are
defined as follows in which we drop the index v. Assume
$\pi\subset Ind_{M_{\theta}(N_{\theta}\cap M)\uparrow
M}\sigma\otimes 1$, where
${\bf{M}}_{\theta}({\bf{N}}_\theta\cap{\bf{M}})$ is a parabolic
subgroup of ${\bf{M}}$ defined by a subset $\theta\subset\Delta$,
the set of simple roots of ${\bf{A}}_0$. Let
$\theta'=\widetilde{w}_0(\theta)\subset\Delta$ and fix a reduced
decomposition $\widetilde{w}_0=\widetilde{w}_{n-1}\cdots
\widetilde{w}_1$ of $\widetilde{w}_0$ (Lemma 2.1.1 of [41]). For
each j, there exists a unique root $\alpha_j\in\Delta$ such that
$\widetilde{w}_j(\alpha_j)<0$. For each j, $2\leq j\leq n-1$, let
$\overline{w}_j=\widetilde{w}_{j-1}\cdots \widetilde{w}_1$. Set
$\overline{w}_1=1.$  Let  $\Omega_j=\theta_j\cup \{\alpha_j\}$,
where $\theta_1=\theta,\ \theta_n=\theta'$, and
$\theta_{j+1}=\widetilde{w}_j(\theta_j),\ 1\leq j\leq n-1$. Then
${\bf{M}}_{\Omega_j}$ contains
${\bf{M}}_{\theta_j}({\bf{N}}_{\theta_j}\cap{\bf{M}}_{\Omega_j})$
as a maximal parabolic subgroup and $\overline{w}_j(\sigma)$ is a
representation of $M_{\theta_j}$. The L-group ${}^LM_{\theta}$
acts on the space of $r_i$, but no longer necessarily
irreducibly. Given an irreducible constituent of this action,
there exists a unique j, $1\leq j\leq n-1$, which under
$\overline{w}_j$ is equivalent to an irreducible constituent of
the action of ${}^LM_{\theta_j}$ on the Lie algebra of the
L-group of ${\bf{N}}_{\theta_j}\cap {\bf{M}_{\Omega_j}}$. Let
$i(j)$ be the index of this subspace and denote by $r'_{i(j)}$
the action of ${}^LM_{\theta_j}$ on it. Finally, let $S_i$ denote
the set of all such j's for a given i.} \rm{(}\it{See Theorem 3.5
and Section 7 of [42]. Also see the discussion just before
Proposition 5.2 of [2.8].$)$ \rm

{\bf Remark 4.5.} If ${\bf{G}}=GL_{t+n},{\bf{M}}=GL_t\times GL_n$
and $\pi=\otimes_v\pi_v$ and $\pi'=\otimes_v\pi'_v$ are cuspidal
representations of $GL_t({\mathbb{A}}_F$ and
$GL_n({\mathbb{A}}_F)$, then $m=1$ and
$L(s,\pi\otimes\widetilde{\pi}',r_1)$ is precisely the
Rankin-Selberg product $L$-function $L(s,\pi\times\pi')$ attached
to ($\pi,\pi'$) (cf. [21,43,44]). In this case each of the local
$L$-functions and root numbers are precisely those of Artin
through parametrization which is now available for $GL_N(F_v)$
for any $N$ due to the work Harris-Taylor [18] and Henniart [19].
As we explain later, this will also be the case for many of our
local factors as a result of our new cases of functoriality which
we shall soon explain. This is quite remarkable, since our
factors are defined using harmonic analysis, as opposed to the
very arithmetic nature of the definition given for Artin factors.
This is a perfect example of how deep Langlands' conjectures are.

{\bf Remark 4.6.} The multiplicativity of local factors, in the
sense of Theorem 3.4, are absolutely crucial in establishing our
new cases of functoriality throughout our proofs [12,23,28]. In
fact, not only do we need them to prove our strong transfers,
they are also absolutely necessary in establishing our weak ones.

\section{Twists by highly ramified characters, holomorphy and boundedness} \label{section 5} \setzero

\vskip-5mm \hspace{5mm}

While the functional equations developed from our method are in
perfect shape and completely general, nothing that general can be
said about the holomorphy and possible poles of these
$L$-functions. On the other hand, there has recently been some
remarkable new progress on the question of holomorphy of these
$L$-function, mainly due to Kim [24,25,31]. They rely on reducing
the existence of the poles to that of existence of certain
$unitary$ automorphic forms, which in turn points to the existence
of certain local unitary representations. One then disposes of
these representations, and therefore the pole, by checking the
corresponding unitary dual of the local group. In view of the
functional equation, this needs to be checked only for
$\mathrm{Re}(s)\geq 1/2$. In fact, to carry this out, one needs to
verify that:
\begin{eqnarray}
& &\mbox{\it Certain local normalized $($as in [41]$)$
intertwining
operators} \nonumber \\
& &\mbox{\it are holomorphic and non-zero for $\mathrm{Re}(s)\geq1/2$}, \label{5.1}
\end{eqnarray}
in each case [24,25,31]. The main issue is that one cannot always get such a contradiction and rule out the pole.
In fact, there are many unitary duals whose complementary series extend all the way to Re($s$)=1.

On the other hand, if one considers a highly ramified twist
$\pi_{\eta}$ (see Theorem 5.1 below) of $\pi$, then it can be
shown quite generally that every $L(s,\pi_{\eta},r_i)$ is entire
(cf. [45] for its local analogue). In fact, if $\eta$ is highly
ramified, then $w_0(\pi_{\eta})\ncong\pi_{\eta}$, whose negation
is a necessary condition for $M(s,\pi_{\eta})$ to have a pole, a
basic fact from Langlands spectral theory of Eisenstein series
(Lemma 7.5 of [33]). This was used by Kim [24], and in view of the
present powerful converse theorems [8,9], that is all one needs to
prove our cases of functoriality [12,23,28,30]. To formalize this,
we borrow the following proposition (Proposition 2.1) from [28],
in order to state the result. It is a consequence of our general
induction (Propositions 4.1 and 4.3) and [24].

{\bf Theorem 5.1.} \it Assume $(5.1)$ is valid. Then there exists
a rational character $\xi\in X({\bf{M}}_F)$ with the following
property: Let S be a non-empty finite set of finite places of F.
For every globally generic cuspidal representation $\pi$ of
${M=\bf{M}}(\mathbb{A}_F)$, there exist non-negative integers
$f_v,\ v\in S$, depending only on the local central characters of
$\pi_v$ for all $v\in S$, such that for every gr\"{o}ssencharacter
$\eta=\otimes_v\eta_v$ of F for which conductor of $\eta_v,\ v\in
S$, is larger than or equal to $f_v$, every L-function
$L(s,\pi_\eta,r_i)$, $i=1,\cdots,m$, is entire, where
$\pi_{\eta}=\pi\otimes(\eta\cdot\xi)$. The rational character
$\xi$ can be simply taken to be $\xi(m)=\mathrm{det}(Ad(m)|{\frak
n}),\ m\in{\bf{M}}$, where $\frak n$ is the Lie algebra of
$\bf{N}$. \rm

The last ingredient in applying converse theorems is that of
boundedness of each $L(s,\pi,r_i)$ in every vertical strip of
finite width, away from its poles, which are finite in number,
again using the functional equation and under Assumption (5.1).
This was proved in full generality by Gelbart-Shahidi [15], using
the theory of Eisenstein series via [33] and [36]. The main
theorem of [15] (Theorem 4.1) is in full generality, allowing
poles for $L$-functions. Here we will state the version which
applies to our $\pi_{\eta}$.

{\bf Theorem 5.2.} \it Under Assumption $(5.1)$, let $\xi$ and
$\eta$ be as in Theorem 5.1. Assume $\eta$ is ramified enough so
that each $L(s,\pi_{\eta},r_i)$ is entire. Then, given a finite
real interval I, each $L(s,\pi_{\eta},r_i)$ remains bounded for
all s with $\mathrm{Re}(s)\in I$}.\rm

The main difficulty in proving Theorem 5.2 is having to deal with
reciprocals of each $L(s,\pi,r_i),\ 2\leq i\leq m$, near and on
the line Re($s$)=1, the edge of the critical strip, whenever
$m\geq2$, which is unfortunately the case for each of our cases of
functoriality. We handle this by appealing to equations (3.5) and
estimating the non-constant term (3.4) by means of [33,36].

\section{New cases of functoriality} \label{section 6}\setzero

\vskip-5mm \hspace{5mm}

Langlands functoriality predicts that every homomorphism between
$L$-groups of two reductive groups over a number field, leads to a
canonical correspondence between automorphic representations of
the two groups. The following instances of functoriality are quite
striking and are consequences of applying recent ingenious
converse theorems of Cogdell and Piatetski-Shapiro [8,9] to
certain classes of $L$-functions whose necessary properties are
obtained mainly from our method. (See [20] for an insightful
survey.) We refer to [11] for more discussion of these results and
the transfer from
$GL_2({\mathbb{A}}_F){\times}GL_2({\mathbb{A}}_F)$ to
$GL_4({\mathbb{A}}_F)$, using Rankin-Selberg method by
Ramakrishnan [37]. (See [23] for a proof using our method.)

{\bf 6.a.} Let $\pi_1=\otimes_v\pi_{1v}$ and
$\pi_2=\otimes_v\pi_{2v}$ be cuspidal representations of
$GL_2({\mathbb{A}}_F)$ and $GL_3({\mathbb{A}}_F)$, respectively.
For each $v$, let $\rho_{iv}$ be the homomorphism of Deligne-Weil
group into $GL_{i+1}(\mathbb{C})$, parametrizing $\pi_{iv},\
i=1,2$. Let $\pi_{1v}\boxtimes \pi_{2v}$ be the irreducible
admissible representation of $GL_6(F_v)$ attached to
$\rho_{1v}\otimes\rho_{2v}$ via [18,19]. Set
$\pi_1\boxtimes\pi_2=\otimes_v(\pi_{1v}\boxtimes\pi_{2v})$, an
irreducible admissible representation of $GL_6({\mathbb{A}}_F)$.
Next, let $\pi=\pi_1,\ \pi_v=\pi_{1v}$ and $\rho_v=\rho_{1v}$. Let
${\rm Sym^3}(\pi_v)$ be the irreducible admissible representation
of $GL_4(F_v)$ attached to ${\rm Sym^3}(\rho_v)$ and set ${\rm
Sym^3}(\pi)=\otimes_v {\rm Sym^3}(\pi_v)$, an irreducible
admissible representation of $GL_4({\mathbb{A}}_F)$. We have:

{\bf Theorem 6.1 [28,30].} \it $a)$ The representations
$\pi_1\boxtimes \pi_2$ and $Sym^3(\pi)$ are automorphic.

$b)$ $Sym^3(\pi)$ is cuspidal, unless $\pi$ is either of dihedral or of tetrahedral type. \rm

In view of [9], one needs to show that
$L(s,(\pi_1\boxtimes\pi_2)\times(\sigma\otimes\eta))$ is
{\bf{nice}} in the sense that it satisfies the contentions of
Theorems 4.4.a, 5.1 and 5.2 for a highly ramified
gr\"{o}ssencharacter $\eta$, where $\sigma$ is a cuspidal
representation of $GL_n({\mathbb{A}}_F),\ n=1,2,3,4$, which is
unramified in every place $v$ where either $\pi_{1v}$ or
$\pi_{2v}$ is ramified. In particular for each $v$, one of
$\pi_{1v},\pi_{2v}$ or $\sigma_v$ is in the principal series. It
then follows from multiplicativity (cf. Theorem 4.4) and the main
results of [43,44], that these $L$-functions are equal to certain
$L$-functions defined from our method. More precisely, we can
take $({\bf{G}},{\bf{M}})$ to be: a) ${\bf{G}}=SL_5,\
{\bf{M}}_D=SL_2\times SL_3$; b) ${\bf{G}}=\mathrm{Spin}(10),\
{\bf{M}}_D=SL_3\times SL_2\times SL_2$; c) ${\bf{G}}=E_6^{sc},\
{\bf{M}}_D=SL_3\times SL_2\times SL_3$; d) ${\bf{G}}=E_7^{sc},\
{\bf{M}}_D=SL_3\times SL_2\times SL_4$, according as $n=1,2,3,4$,
respectively. This leads to a proof that $\pi_1\boxtimes\pi_2$ is
weakly automorphic. The strong transfer requires a lot more work,
involving base change, both normal [2] and non-normal [22], and
finally a local result [5]. Automorphy of ${\rm Sym^3}(\pi)$ is a
consequence of applying the first part to ($\pi$,Ad($\pi$)),
where Ad($\pi$) is the adjoint of $\pi$, established by
Gelbart-Jacquet [14]. It does not require the use of [5].

Observe that we have in fact proved that the homomorphisms
$GL_2(\mathbb{C})\otimes GL_3(\mathbb{C})\subset GL_6(\mathbb{C})$
and ${\rm Sym^3}$: $GL_2(\mathbb{C})\rightarrow GL_4(\mathbb{C})$
are functorial. Neither are endoscopic.

{\bf 6.b.} Let $\Pi=\otimes_v\Pi_v$ be a cuspidal representation of $GL_4({\mathbb{A}}_F)$ and let
$\Lambda^2:GL_4(\mathbb{C})\rightarrow GL_6({\mathbb{C}})$ be the exterior square map. Also with $\pi$ as in 6.a,
let ${\rm Sym^4} (\pi)=\otimes_v {\rm Sym^4}(\pi_v)$, where ${\rm Sym^4}(\pi_v)$ is attached to ${\rm
Sym^4}(\rho_v)$. Then

{\bf Theorem 6.2 (cf. [23]).} \it $a)$ The map $\Lambda^2$ is weakly functorial, in the sense that there exists an
automorphic form on $GL_6(\mathbb{A}_F)$ whose local components are equal to $\Lambda^2(\Pi_v)$ for all v, except
if $v|2$ or $v|3$. Here $\Lambda^2(\Pi_v)$ is defined by the local Langlands conjecture [18,19].

$b)$ $Sym^4(\pi)$ is an automorphic representation of $GL_5({\mathbb{A}}_F)$. \rm

We point out that b) is obtained by applying a) to
$\mathrm{Sym}^3(\pi)$. a) is proved by applying our method to Spin
groups (Case $D_n-1$ of [40], $n=k+4,\ k=0,1,2,3$).

{\bf Proposition 6.3 (cf. [29]).} \it $Sym^4(\pi)$ is cuspidal
unless $\pi$ is either of dihedral, tetrahedral or octahedral
type. \rm

Let $\pi=\otimes_v\pi_v$ be a cusp form on $GL_2(\mathbb{A}_F)$.
For each unramified $v$, let $\alpha_v$ and $\beta_v$ be the Hecke
eigenvalues of $\pi_v$. Then as corollary to Proposition 6.3 we
have the following striking improvements towards Ramanujan and
Selberg conjectures.

{\bf Corollary 6.4.} \it $a)$ $(${cf. [29]}$)$  Assume F is an arbitrary number field. Then
$q_v^{-1/9}<|\alpha_v|$ and $|\beta_v|<q_v^{1/9}$. b) (cf. [27]). Assume $F=\mathbb{Q}$. Then
$p^{-7/64}\leq|\alpha_p|$ and $|\beta_p|\leq p^{7/64}$. Similar estimates are valid for the Selberg conjecture.
More precisely, the smallest positive eigenvalue $\lambda_1(\Gamma)$ of the Laplace operator on $L^2(\Gamma\
{\frak n})$ for every congruence subgroup $\Gamma$ satisfies
$\lambda_1(\Gamma)\geq\dfrac{975}{4096}\cong0.2380\cdots$ \rm

{\bf 6.c.} Let $i:Sp_{2n}(\mathbb{C})\hookrightarrow GL_{2n}(\mathbb{C})$ be the natural embedding. Let
$\pi=\otimes_v\pi_v$ be a generic cuspidal representation of $SO_{2n+1}({\mathbb{A}}_F)$. For each unramified $v$,
let $\{A_v\}\subset Sp_{2n}(\mathbb{C})$ be the Hecke-Frobenius conjugacy class parametrizing $\pi_v$. Let
$i(\pi_v)$ be the unramified representation of $GL_{2n}(F_v)$ attached to $\{i(A_v)\}$. Then the main theorem of
[12] proves:

{\bf Theorem 6.5 [12].} \it The embedding i is weakly functorial,
i.e., there exist an automorphic representation of
$GL_{2n}(\mathbb{A}_F)$ whose components are equal to $i(\pi_v)$
for almost all $v$. \rm

This is proved by applying our method to maximal parabolics of
appropriate odd special orthogonal groups (Case $B_n$ of [40]).
The strong transfer is now also established by
Ginzburg-Rallis-Soudry [16] as well as Kim [26] by building upon
Theorem 6.5.

{\bf Final Comments.} Many other cases are in progress. Among
them are a proof of the existence of Asai transfer [32] using our
method, which was originally proved by Ramakrishnan [38], using
the Rankin-Selberg method. This is the first case where one needs
to use quasisplit groups. Since the issue of stability of root
numbers [10] (cf. [11]) seems to be close to being settled by
means of our method [46], many others transfers should now be
available. A similar approach for nongeneric representations was
initiated in [13].

\label{lastpage}


\begin{thebibliography}{aa}
\bibitem{1} J. Arthur, {\it Endoscopic $L$-functions and a combinatorial indentity}, Dedicated to H.S.M.
Coxeter, Canad. J. Math., {\bf{51}} (1999), 1135--1148.

\bibitem{2} J. Arthur and L. Clozel, {\it Simple algebras,
Base change, and the Advanced Theory of the Trace Formula},
Annals of Math. Studies, no. 120, Princeton University Press,
1989.

\bibitem{3} M. Asgari, {\it Local $L$-functions for split
spinor groups}, Canad. J. Math., (to appear).

\bibitem{4} A. Borel., {\it Automorphic L-functions},
Automorphic Forms and Automorphic Representations, Proc. Sympos.
Pure Math., vol. 33; II, Amer. Math. Soc., Providence, RI, 1979,
 27--61.

\bibitem{5} C. J. Bushnell and G. Henniart, {\it On certain
dyadic representations}, Annals of Math., Appendix to [28], (to
appear).

\bibitem{6} W. Casselman and F. Shahidi, {\it On
irreducibility of standard modules for generic representations},
Ann. Scient. \'{E}c.Norm.Sup., {\bf{31}} (1998), 561--589.

\bibitem{7} W. Casselman and J. A. Shalika, {\it The
unramified principal series of p-adic groups II; The Whittaker
function}, Comp. Math., {\bf{41}} (1980), 207--231.

\bibitem{8} J. W. Cogdell and I. I. Piatetski-Shapiro,
{\it Converse theorems for $GL_n$}, Publ. Math. IHES, {\bf{79}}
(1994), 157--214.

\bibitem{9} \underline{\hskip10mm}, {\it Converse theorems
for $GL_n$ II}, J. Reine Angew. Math., {\bf{507}} (1999),
165--188.

\bibitem{10} \underline{\hskip10mm}, {\it Stability of gamma
factors for $SO(2n+1)$}, Manuscripta Math., {\bf{95}} (1998),
437--461.

\bibitem{11} \underline{\hskip10mm}, {\it Converse Theorems,
Functoriality and Applications to Number Theory}, These Proceedings.

\bibitem{12} J. W. Cogdell, H. Kim, I. I. Piatetski-Shapiro and
F. Shahidi, {\it On lifting from classical groups to $GL_N$},
Publ. Math. IHES, {\bf{93}} (2001), 5--30.

\bibitem{13} S. Friedberg and D. Goldberg, {\it On local
coefficients for nongeneric representations of some classical
groups}, Comp. Math., {\bf{116}} (1999), 133--166.

\bibitem{14} S. Gelbart and H. Jacquet, {\it A relation
between automorphic representations of $GL(2)$ and $GL(3)$}, Ann.
Scient. \'{E}c. Norm. Sup., {\bf{11}}( 1978), 471--552.

\bibitem{15} S. Gelbart and F. Shahidi, {\it Boundedness of
automorphic L-functions in vertical strips}, {Journal of AMS},
{\bf 14} (2001), 79--107.

\bibitem{16} D. Ginzburg, S. Rallis and D. Soudry, {\it Generic
automorphic forms on SO(2n+1): functorial lift to $GL(2n)$,
endoscopy and base change}, IMRN {\bf{14}} (2001), 729--764.
\bibitem{17} Harish-Chandra, {\it Automorphic forms on
semisimple Lie groups}, {SLN {\bf{62}} (1968), Berlin-Heidelberg-New York.}

\bibitem{18} M. Harris and R. Taylor, {\it On the geometry
and cohomology of some simple Shimura varieties}, Annals of Math.
Studies, no. 151, Princeton University Press, 2001.

\bibitem{19} G. Henniart, {\it Une preuve simple des
conjectures de Langlands pour GL(n) sur un corps p-adique},
Invefnt. Math., {\bf{139}} (2000), 439--455.

\bibitem{20} G. Henniart, {\it Progr\`{e}s r\'{e}cents en
fonctorialit\'{e} de Langlands}, Seminaire Bourbaki, Juin 2001, Expos\'{e} 890, 890-1 to 890-21.

\bibitem{21} H. Jacquet, I. Piatetski-Shapiro and J. Shalika,
{\it Rankin-selberg convolutions}, Amer. J. Math., {\bf{105}}
(1983), 367--464.

\bibitem{22} \underline{\hskip10mm}, {\it Relvement cubique
non normal}, C. R. Acad. Sci. Paris Sr. I Math., {\bf 292}, no.
12 (1981), 567--571.

\bibitem{23} H. Kim, {\it Functoriality for the exterior
square of $GL_4$ and symmetric fourth of $GL_2$}, preprint (2000).

\bibitem{24} \underline{\hskip10mm}, {\it Langlands-Shahidi
method and poles of automorphic L-functions: Application to
exterior square L-functions}, Can. J. Math., {\bf{51}} (1999),
835--849.

\bibitem{25} \underline{\hskip10mm}, {\it Langlands-Shahidi
method and poles of automorphic L-functions II}, Israel J. Math.,
{\bf{117}} (2000), 261--284.

\bibitem{26} \underline{\hskip10mm}, {\it Residual spectrum
of odd orthogonal groups}, IMRN, {\bf{17}} (2000), 873--906.

\bibitem{27} H. Kim and P. Sarnak, {\it Refined estimates
towards the Ramanujan and Selberg conjectures}, Appendix 2 to [23].

\bibitem{28} H. Kim and F. Shahidi, {\it Functorial products
for $GL_2\times GL_3$ and the symmetric cube for $GL_2$}, Annals
of Math., (to appear).

\bibitem{29} \underline{\hskip10mm}, {\it Cuspidality of
symmetric powers with applications}, Duke Math. J., {\bf{112}}
(2002), 177--197.

\bibitem{30} \underline{\hskip10mm}, {\it Functorial
products for $GL_2\times GL_3$ and functorial symmetric cube for
$GL_2$}, C.R. Acad. Sci. Paris, {\bf{331}} (2000), 599--604.

\bibitem{31} \underline{\hskip10mm}, {\it Symmetric cube
L-functions for $GL_2$ are entire}, Ann. of Math., {\bf{150}}
(1999), 645--662.

\bibitem{32} M. Krishnamurthy, The weak Asai transfer to
$GL$(4) via Langlands-Shahidi method, Thesis, Purdue University
(2002).

\bibitem{33} R. P. Langlands, {\it On the Functional
Equations Satisfied by Eisenstein Series}, Lecture Notes in Math., Vol 544, Springer-Verlag, 1976.

\bibitem{34} \underline{\hskip10mm}, {\it Euler
Products}, Yale University Press, 1971.

\bibitem{35} C. Moeglin and J.-L. Waldspurger, {\it Spectral
decomposition and Eisenstein series}, Cambridge Tracts in Math., vol. 113, Cambridge University Press, 1995.

\bibitem{36} W. M\"{u}ller, {\it The trace class conjecture
in the theory of automorphic forms}, Ann. of Math., {\bf{130}}
(1989), 473--529.

\bibitem{37} D. Ramakrishnan, {\it Modularity of the
Rankin-Selberg L-series, and multiplicity one for SL}(2), Ann. of
Math., {\bf{152}} (2000), 45--111.

\bibitem{38} \underline{\hskip10mm}, {\it Modularity of
solvable Artin representations of GO}(4)-$type$, IMRN, {\bf{1}}
(2002), 1--54.

\bibitem{39} F. Shahidi, {\it Functional equation satisfied
by certain L-functions}, Comp. Math., {\bf{37}}(1978), 171--207.

\bibitem{40} \underline{\hskip10mm}, {\it On the Ramanujan
conjecture and finiteness of poles for certain L-functions}, Ann.
of Math., {\bf{127}} (1988), 547--584.

\bibitem{41} \underline{\hskip10mm}, {\it On certain
L-functions}, Amer, J. Math., {\bf{103}}(1981), 297--355.

\bibitem{42} \underline{\hskip10mm}, {\it A  proof of
Langlands conjecture on Plancherel measures;  Comple-\\
mentary series for p-adic groups}, Annals of Math., {\bf{132}}
(1990), 273--330.

\bibitem{43} \underline{\hskip10mm}, {\it Local coefficients
as Artin factors for real groups}, Duke Math. J., {\bf{52}}
(1985), 973--1007.

\bibitem{44} \underline{\hskip10mm}, {\it Fourier transforms
of intertwining operators and Plancherel measures for $GL(n)$},
Amer. J. of Math., {\bf{106}} (1984), 67--111.

\bibitem{45} \underline{\hskip10mm}, {\it Twists of a
general class of L-functions by highly ramified characters},
Canad. Math. Bull., {\bf{43}} (2000), 380--384.

\bibitem{46} \underline{\hskip10mm}, {\it Local coefficients
as Mellin transforms of Bessel functions; Towards a general
stability}, preprint (2002).
\end{thebibliography}
\end{document}